\documentclass{amsart}

\usepackage[frame,cmtip,arrow,matrix,line,graph,curve]{xy}
\usepackage{amsmath,amssymb, amscd}
\usepackage{latexsym}
\usepackage{graphics}

\numberwithin{equation}{section}

\newtheorem{theorem}{Theorem}[section]
\newtheorem{proposition}[theorem]{Proposition}

\newtheorem{lemma}[theorem]{Lemma}

\newtheorem{remit}[theorem]{Remark}
\newtheorem{definit}[theorem]{Definition}

\newenvironment{definition}{\begin{definit}\rm}{\end{definit}}

\newenvironment{remark}{\begin{remit}\rm}{\end{remit}}

\newcommand{\pp}{\mathbb{P}}

\newcommand{\cc}{\mathbb{C}}

\newcommand{\zz}{\mathbb{Z}}

\newcommand{\git}{/\!\!/}

\newcommand{\cA}{\mathcal{A} }
\newcommand{\cB}{\mathcal{B} }
\newcommand{\cF}{\mathcal{F} }
\newcommand{\cE}{\mathcal{E} }
\newcommand{\cM}{\mathcal{M} }
\newcommand{\cN}{\mathcal{N} }
\newcommand{\cL}{\mathcal{L} }
\newcommand{\cO}{\mathcal{O} }

\newcommand{\cU}{\mathcal{U} }
\newcommand{\End}{\mathrm{End} }

\newcommand{\tcN}{\tilde{\cN} }

\newcommand{\tcF}{\tilde{\cF} }
\newcommand{\cD}{\mathcal{D} }

\newcommand{\tcD}{\tilde{\cD} }
\newcommand{\tD}{\tilde{D} }

\newcommand{\tsig}{\tilde{\Sigma} }

\newcommand{\tdel}{\tilde{\Delta} }

\newcommand{\tcL}{\tilde{\mathcal{L}} }

\newcommand{\e}{\epsilon }
\newcommand{\s}{\sigma }

\newcommand{\fK}{\mathfrak{K} }
\newcommand{\fR}{\mathfrak{R} }
\newcommand{\bE}{\mathbf{E} }
\newcommand{\bK}{\mathbf{K} }
\newcommand{\bN}{\mathbf{N} }
\newcommand{\bS}{\mathbf{S} }
\newcommand{\cDoo}{\mathcal{D}^{(1)}_1 }
\newcommand{\cDto}{\mathcal{D}^{(1)}_2 }
\newcommand{\cDtt}{\mathcal{D}^{(2)}_2 }

\begin{document}

\title[Moduli of Hecke cycles]{Cohomology of the moduli space of Hecke cycles}
\date{}

\author{Insong Choe, Jaeyoo Choy and Young-Hoon Kiem}
\address{School of Mathematics, KIAS, 207-43 Cheongnyangni 2-dong, 130-722 Korea}
\email{ischoe@kias.re.kr}
\address{Dept of Mathematics, Seoul National University,
Seoul, 151-747 Korea} \email{kiem@math.snu.ac.kr}
\email{donvosco@math.snu.ac.kr}

\thanks{Young-Hoon Kiem was partially supported by KOSEF}


\begin{abstract}
Let $X$ be a smooth projective curve of genus $g \ge 3$ and let
$M_0$ be the moduli space of semistable bundles over $X$ of rank
$2$ with trivial determinant. Three different desingularizations
of $M_0$ have been constructed by Seshadri \cite{Se1},
Narasimhan-Ramanan \cite{NR}, and Kirwan \cite{k5}. In this paper,
we construct a birational morphism from Kirwan's desingularization
to Narasimhan-Ramanan's, and prove that the Narasimhan-Ramanan's
desingularization (called the moduli space of Hecke cycles) is the
intermediate variety between Kirwan's and Seshadri's as was
conjectured recently in \cite{KL}. As a by-product, we compute the
cohomology of the moduli space of Hecke cycles.
\end{abstract}
 \maketitle


\section{Introduction}

Let $X$ be a smooth projective curve of genus $g \ge 3$ over the
complex number field. Let $M_0$ be the moduli space of semistable
bundles over $X$ of rank $2$ with trivial determinant. Then $M_0$
is a singular normal projective variety of dimension $3g-3$. Its
singular locus is the Kummer variety $\fK$ which consists of the
S-equivalence classes of strictly semistable bundles $E = L \oplus
L^{-1}$ for $L \in Pic^0 (X)$.

There are three different constructions to desingularize $M_0$:

(1) Seshadri's desingularization $\bS$ (\cite{Se1}),

(2) Narasimhan-Ramanan's desingularization $\bN$ (\cite{NR}),
called the \emph{moduli space of Hecke cycles}, and

(3) Kirwan's desingularization $\bK$ (\cite{k5}).

The first two desinglarizations $\bS$ and $\bN$ come from certain
moduli problems, while $\bK$ is obtained as a result of more
general construction of a partial desingularization of a GIT
quotient, which was studied by F. Kirwan in \cite{k2}.

Recently Y.-H. Kiem and J. Li in \cite{KL} constructed a morphism
$f: \bK \rightarrow \bS$ and described it explicitly as a
composition of two blow-downs:
$$
\xymatrix{
f:&\bK\ar[r]^{f_{\sigma}}&\bK_{\sigma}\ar[r]^{f_{\epsilon}}&\bK_{\e}}
(\cong \bS).
$$
Also they conjectured that the intermediate variety $\bK_{\sigma}$
is isomorphic to the moduli space of Hecke cycles $\bN$
(\cite{KL}, Conjecture 5.7). In this paper, we give a proof of
this conjecture and compute the cohomology of $\bN$ as its
by-product. For this, we construct a birational morphism (Theorem
\ref{mainthm})
$$
\rho : \bK \rightarrow \bN
$$
and then show that this coincides with the morphism $ f_\sigma :
\bK \rightarrow \bK_\sigma $ of \cite{KL} by examining the fibers
of $\rho$ (Proposition \ref{iso7}). M.S. Narasimhan and S. Ramanan
conjectured that the desingularization $\bN$ can be blown down
along certain projective fibrations to obtain another nonsingular
model of $M_0$ (\cite{NR}, page 292) and this was proved by N.
Nitsure \cite{Ni}. Our result shows that this blown down process
corresponds to the morphism
$$
f_\epsilon : \ \bK_\sigma (\cong \bN) \longrightarrow \bK_\epsilon
(\cong \bS).
$$
In summary, the three desingularizations are related by morphisms
$$\bK\to
\bN\to \bS$$ which can be described explicitly as blow-up maps
along smooth subvarieties.

The strategy of the construction of $\rho$ is similar as that of
$f$ in \cite{KL}. There is a birational map $\rho' : \bK
\dashrightarrow \bN$ which is defined on the open subset $M_0^s$
of stable bundles. By GAGA and Riemann's extension theorem
\cite{Kran}, it suffices to show that $\rho'$ can be extended to a
continuous map with respect to the usual complex topology. By
Luna's slice theorem, for each point $x \in M_0 \backslash M_0^s$,
there is an analytic submanifold $W$ of the Quot scheme whose
quotient by the stabilizer $H$ of a point in both $W$ and the
closed orbit represented by $x$ is analytically equivalent to a
neighborhood of $x$ in $M_0$. Furthermore, Kirwan's
desingularization $\tilde{W} \git H$ of $W \git H$ is a
neighborhood of the preimage of $x$ in $\bK$.

There is a universal family $\cU$ of rank 2 vector bundles over
$X$ parameterized by $\tilde{W}$, which is induced from the
universal bundle over the Quot scheme. By applying an elementary
modification with respect to the points of the curve $X$, we have
a family $\cU'$ of rank 2 vector bundles of determinant $\cO_X
(-x)$ for some $x \in X$, which is parameterized by the projective
bundle $\pp \cU^*$ over $\tilde{W} \times X$. For any point $w \in
\tilde{W}$ lying over a stable bundle in $M_0$, the bundles of
$\cU'$ parameterized by the fiber of $w$ are all stable, and a
good Hecke cycle is associated to $w$. This process yields the
birational map $ \rho' : \bK \dashrightarrow \bN$.

The problem is that for the points $w \in \tilde{W}$ lying over a
strictly semistable bundle in $M_0$, some points of $\pp \cU^*$ in
the fiber of $w$ parameterize unstable bundles in $\cU'$. To
remedy this, we blow up $\pp \cU^*$ and then apply an elementary
modification of $\cU'$ along the exceptional divisors. Local
computations of the transition data show that the resulting family
$\cU''$ yields an analytic extension $ \rho : \bK \rightarrow \bN$
of $\rho'$.

This paper is organized as follows. In section 2, we explain the
elementary modification of vector bundles, focusing on its local
computations which will be used repeatedly in this paper. In
section 3 and section 4, we briefly review the
Narasimhan-Ramanan's and Kirwan's desingularizations respectively.
In section 5 and section 6, we construct the birational morphism
$\rho : \bK \rightarrow \bN$. In section 7, we examine the fibers
of $\rho$ and prove that $\rho$ is in fact a blow-up along a
smooth subvariety of $\bN$. In section 8, we compute the
cohomology of $\bN$ using the morphism $\rho$. We remark that N.
Nitsure(\cite{Ni}) computed the third cohomology group $H^3(\bN,
\zz)$ of $\bN$.

\section{Elementary modification}

Let $X$ be a smooth projective curve over the complex number
field. Let $E$ be a vector bundle over $X$ and $E_x$ the fiber of
$E$ at $x$. For simplicity, assume ${\rm rk} (E) = 2$.

For any nonzero homomorphism $\nu: E_x \rightarrow \cc$, we have
an exact sequence
\begin{equation} \label{Hecke1}
0 \rightarrow E^{\nu} \rightarrow E \stackrel{\nu}{\rightarrow}
 \cc_x \rightarrow 0,
\end{equation}
where $\cc_x$ is the skyscraper sheaf supported at $x$. Then
$E^{\nu} = \ker(\nu)$ is locally free and is called an {\it
elementary modification} of $E$.

In terms of the transition matrices, this process can be described
as follows. Choose a local trivialization of $E$ with an open
covering $\{ V_i \}$ of $X$ and the transition matrices
\begin{equation}
\{ \  g_{ij} = \left(\begin{matrix} a_{ij} & b_{ij}\\
c_{ij}& d_{ij} \end{matrix}\right) : V_i \cap V_j \
\longrightarrow \ GL(2, \cc) \ \}.
\end{equation}
We can refine the covering so that $x$ is contained in $V_1$ only.
Let $\zeta$ be a coordinate function on $V_1$ such that
$\zeta(x)=0$.

Suppose that $\nu: E_x \cong \cc^2 \rightarrow \cc$ is the first
projection. Then a local section $(f,g)$ of the sheaf $E^{\nu}$ on
$V_1$ is $(\zeta f, g)$ when considered as a local section of $E$
on $V_1$. Hence from the computation
\begin{equation}
\left(\begin{matrix} f\\ g\end{matrix}\right)\leftrightarrow
\left(\begin{matrix} \zeta f\\ g\end{matrix}\right) \mapsto
g_{1j} \left(\begin{matrix} \zeta f\\
g\end{matrix}\right) = \left(\begin{matrix} \zeta a_{1j} f + b_{1j} g\\
\zeta c_{1j} f + d_{1j} g\end{matrix}\right),
\end{equation}
the transition matrix of $E^{\nu}$ from $V_1$ to $V_j$ for $j\ne
1$ is
$$
\left(\begin{matrix}
\zeta a_{1j} & b_{1j}\\
\zeta c_{1j} & d_{1j}\end{matrix}\right).
$$
Also, the transition of $E^{\nu}$ from $V_j$ to $V_1$ for $j\ne 1$
is the inverse matrix
$$
\left(\begin{matrix}
\zeta^{-1} a_{j1} & \zeta^{-1} b_{j1}\\
 c_{j1} & d_{j1} \end{matrix}\right)
$$
and the other transition matrices are unchanged. Note that
$E^{\nu} \cong E^{\lambda \nu}$ for any nonzero $\lambda \in \cc$.

In this way, we can produce vector bundles $E^{\nu}$ of
determinant $L(-x)$ from $E$ of determinant $L$. In \cite{NR},
Narasimhan and Ramanan used this process to construct the Hecke
cycles, as will be reviewed in next section.

Later we will also use the elementary modification to construct a
morphism $\rho:\bK \rightarrow \bN$. It requires the following
generalization to higher dimensions. Let $S$ be a smooth complex
manifold and let $Z$ be a smooth hypersurface of $S$. Let $E$
(resp. $F$) be a vector bundle on $S$ (resp. $Z$) with ${\rm
rk}(F) < {\rm rk} (E)$. Assume that there is a surjective
homomorphism $\nu : E|_Z \rightarrow F$. Then the kernel $E^{\nu}$
of the composition $E \rightarrow E|_Z \xrightarrow{\nu} F$ is
locally free and defines a vector bundles on $S$. This situation
can be summarized in the following diagram(see \cite{Ma}).
\begin{equation}
\begin{CD}
@. 0 @. 0 @. @. \\
@. @AAA @AAA @. @. \\
0 @>>> \ker (\nu) @>>> E|_Z @>\nu>>  F @>>> 0 \\
@.     @AAA         @AAA             @|   @. \\
0 @>>> E^{\nu}   @>>>  E  @>>> F @>>> 0 \\
@. @AAA @AAA @. @. \\
@. E(-Z) @=  E(-Z) @. @. \\
@. @AAA @AAA @. @. \\
@. 0 @. 0 @. @. \\
\end{CD}
\end{equation}
\vskip 10pt

Now let $X$ be an algebraic curve as before, $S$ a complex
manifold, and $E \rightarrow S \times X$ a family of vector
bundles over $X$ parameterized by $S$. For simplicity, assume
$\dim S = 2$. Let $\pi: \tilde{S} \rightarrow S$ be the blow-up at
one point $\theta \in S$ with the exceptional divisor $Z$. Suppose
that $E_{\theta} \cong L_1 \oplus L_2$ for some line bundles $L_1$
and $L_2$ on $X$. Let $\tilde{E} := (\pi \times 1_X)^* E$ and
$\tilde{L}_i := (\pi \times 1_X)^* L_i$ ($i=1,2$) be the families
of bundles parameterized by $\tilde{S}$ and $Z$ respectively so
that $\tilde{E}|_{Z \times X} \cong \tilde{L}_1 \oplus
\tilde{L}_2$. Consider $\tilde{E}$ over $\tilde{S} \times X$(resp.
$\tilde{L}$ over $Z \times X$) as playing the role of $E$ over
$S$(resp. $F$ over $Z$) in the above. Then we have
\begin{lemma}
Let $\nu$ be the first (resp. second) projection $\tilde{E}|_{Z
\times X} \rightarrow \tilde{L}_1$. Then the associated elementary
modification $\tilde{E}^{\nu}$ defines a family of vector bundles
over $X$ such that for each $\tilde{\theta} \in Z$,
$\tilde{E}^{\nu}|_{\tilde{\theta} \times X}$ is an extension of
$L_2$ by $L_1$ (resp. $L_1$ by $L_2$).
\end{lemma}
\begin{proof}
Choose a local coordinate $(z,t)$ of $S$ in a small neighborhood
$U$ of $\theta = (0,0)$. Let $\tilde{\theta} \in Z$ represent the
line $l_\tau : t=\tau z $ in $U$ for some $\tau \in \cc$. Choose
an open covering $\{ V_i \}$ of $X$ such that $E|_{l_\tau \times
V_i}$ are all trivial. Fix a trivialization for each $V_i$ and let
$L_k^\tau = \tilde{L}_k|_{l_\tau \times X}$ for $k=1,2$. Since
$E|_{0 \times X} \cong L_1 \oplus L_2$, the transition matrices of
$\tilde{E}|_{l_\tau \times X}$ are of the form
\begin{equation} \label{un}
\left(\begin{matrix}
\lambda_{ij} & zb_{ij}\\
z c_{ij} & \mu_{ij}\end{matrix}\right)
\end{equation}
where $\{ \lambda_{ij}|_{z=0} \}$ and $\{ \mu_{ij}|_{z=0} \}$ are
the transition functions of $L_1$ and $L_2$ respectively. From the
properties of transition maps, we have
$$
(\lambda_{jk} \mu_{jk}^{-1}) (\mu_{ij}^{-1} b_{ij}) +
(\mu_{jk}^{-1} b_{jk}) = (\mu_{ki}^{-1} b_{ki}).
$$
This shows that the data $\{ \mu_{ij}^{-1} b_{ij} |_{z=0} \}$
define a C\v{e}ch cocycle in $H^1 (X, L_1 \otimes L_2^{-1})$.
Similarly, the data $\{ \lambda_{ij}^{-1} c_{ij} |_{z=0} \}$
define a C\v{e}ch cocycle in $H^1 (X, L_1^{-1} \otimes L_2)$.

The modified bundle $\tilde{E}^\nu$ over $l_\tau \times X$ is
given by the kernel of the composition
$$
\tilde{E}|_{l_\tau \times X} \cong L_1^\tau \oplus L_2^\tau
\rightarrow L_1^\tau.
$$
Note that any section of $\tilde{E}^\nu$ over $l_\tau \times V_i$
is of the form $(zf, g)$ when considered as a section of
$\tilde{E}$. From the computation
$$
\left(\begin{matrix} f\\ g\end{matrix}\right)\leftrightarrow
\left(\begin{matrix}z f\\ g\end{matrix}\right)\mapsto
\left(\begin{matrix}
\lambda_{ij} & zb_{ij} \\
z c_{ij} & \mu_{ij} \end{matrix}\right) \left(\begin{matrix}z f\\
g\end{matrix}\right)
= \left(\begin{matrix} z(\lambda_{ij}f+b_{ij} g)\\
z^2c_{ij}f+ \mu_{ij}g\end{matrix}\right)\leftrightarrow
\left(\begin{matrix}
\lambda_{ij} & b_{ij}\\
z^2c_{ij} & \mu_{ij}\end{matrix}\right)
\left(\begin{matrix} f\\
g\end{matrix}\right)
$$
the transition for $\tilde{E}^\nu |_{\tilde{\theta} \times X}$ is
$$
\left(\begin{matrix}
\lambda_{ij} & b_{ij}\\
0 & \mu_{ij}\end{matrix}\right).
$$
Hence $\tilde{E}^\nu|_{\tilde{\theta} \times X}$ is an extension
of $L_2$ by $L_1$ whose extension class is given by $\{
\mu_{ij}^{-1} b_{ij} |_{z=0} \}$. The same argument proves the
case of the second projection.
\end{proof}


\section{Moduli space of Hecke cycles}

Let $X$ be a smooth projective curve of genus $g \ge 3$ over the
complex number field. Let $M_0$ be the moduli space of semistable
bundles over $X$ of rank $2$ with trivial determinant. Then $M_0$
is a singular normal projective variety of dimension $3g-3$. Its
singular locus is the Kummer variety $\fK$ which consists of the
S-equivalence classes of non-stable bundles $E = L \oplus L^{-1}$
for $L \in Pic^0 (X)$. In \cite{NR}, Narasimhan and Ramanan
constructed a desingularization
$$
\varphi_\bN : \bN \rightarrow M_0
$$
which is an isomorphism over the open subset $M_0^s$ of stable
bundles. The smooth variety $\bN$ is called {\it the moduli space
of Hecke cycles}. In this section, we review its construction.

For any point $x \in X$, let $M_x$ be the moduli space of stable
vector bundles over $X$ of rank $2$ whose determinants are
isomorphic to $\cO_X(-x)$. Let $M_X$ denote the moduli space of
stable bundles over $X$ of rank $2$ whose determinants are
isomorphic to $\cO_X(-x)$ for some $x \in X$, i.e.,
${\displaystyle M_X = \bigcup_{x \in X} M_x}$ inside the moduli
space of stable bundles over $X$ of rank $2$ and degree $1$.

For any stable bundle $E \in M_0^s$ and any $ \nu \in \pp E_x^*$,
we get an associated elementary modification
\begin{equation*}
0 \rightarrow E^\nu \rightarrow E \xrightarrow{\nu} \cc_x
\rightarrow 0
\end{equation*}
so that $\det(E^\nu) = \cO_X(-x)$. Since $E^\nu$ is again
stable(\cite{NR}, Lemma 5.5), $E^\nu \in M_x$. Hence by the
universal property of $M_X$, we have a morphism
$$
\theta_E : \pp E^* \rightarrow M_X.
$$
More generally, for any family $W \rightarrow S \times X$ of
stable bundles in $M_0^s$, there is a canonical morphism $\theta_W
: \pp W^* \rightarrow M_X$. Moreover, this is a closed immersion,
provided that $W_{s_1} \ncong W_{s_2}$ whenever $s_1 \neq s_2 \in
S$ (\cite{NR}, Lemma 5.9).

From this, we get a morphism
$$
\Phi : M_0^s \rightarrow Hilb(M_X)
$$
into the Hilbert scheme of $M_X$, defined by $\Phi(E) = \theta_E
(\pp E^*) \subset M_X$.

\begin{definition} (\cite{NR}, Definition
5.12) \label{goodHecke} For a stable bundle $E \in M_0^s$, the
cycle $\Phi(E)$ in $M_X$ is called the {\it good Hecke cycle
associated to $E$}. Any subscheme in the irreducible component of
$Hilb(M_X)$ containing the good Hecke cycles is called a {\it
Hecke cycle}.
\end{definition}

\begin{theorem}
{\rm (\cite{NR}, Theorem 5.13)} Via the morphism $\Phi$, $M_0^s$
is isomorphic to an open subscheme of $Hilb(M_X)$ consisting of
the good Hecke cycles. $\Box$
\end{theorem}

To compute the Hilbert polynomial of the good Hecke cycles, we fix
an ample line bundle on $M_X$. Let $K_{det}$ denote the canonical
line bundle along the fibers of the fibration $det: M_X
\rightarrow X$. Then $\cO(1) := K_{det}^* \otimes (det)^*K_X$ is
an ample line bundle on $M_X$ (\cite{NR}, Lemma 7.1).

\begin{lemma}
{\rm (\cite{NR}, Lemma 7.2)} The Hilbert polynomial of a good
Hecke cycle is $P(n) = (4n+1)(4n-1)(g-1)$ with respect to
$\cO(1)$. $\Box$
\end{lemma}
Recall that the canonical line bundle of $M_x$ is isomorphic to
$\cL_x^{\otimes (-2)}$ for the ample generator $\cL_x$ of
$Pic(M_x) \cong \zz$ (\cite{Ra}). Also, it is known that $\cL_x$
is very ample (\cite{BV}) and we can think $M_x$ as a projective
variety embedded in $|\cL_x|^*$. In this setting, we see that a
good Hecke cycle in $M_X$ restricts to a conic on $M_x$ for each
$x \in X$.

\begin{theorem}
{\rm (\cite{NR}, \S 8)} Let $\bN$ be the irreducible component of
$Hilb^{P(n)}(M_X)$ containing good Hecke cycles. Then $\bN$ is a
nonsingular variety of dimension $3g-3$. Moreover, there is a
morphism
$$
\varphi_\bN : \bN \rightarrow M_0
$$
which is an isomorphism over the set $M_0^s$ of stable points.
$\Box$
\end{theorem}

The fibers of $\varphi_\bN$ over the boundary locus $M_0
\backslash M_0^s = \fK $ are described as follows (see \cite{NR},
Proposition 7.8 and Theorem 8.14). First consider $L \in Pic^0(X)$
with $L^2 \ncong \cO_X$ and let $l = [L \oplus L^{-1}] \in \fK -
\fK_0$ be a non-nodal point in the Kummer variety in $M_0$. The
fiber $\varphi^{-1}(l)$ is isomorphic to the product of two
$(g-2)$-dimensional projective spaces $\pp H^1(X, L^2)$ and $\pp
H^1 (X, L^{-2})$. Any choice of two points from $\pp H^1(X, L^2)$
and $\pp H^1 (X, L^{-2})$ gives rise to two lines in $\pp H^1(X,
L^2(-x))$ and $\pp H^1 (X, L^{-2}(-x))$. It can be shown that
these two lines meet at the unique intersection point $\pp H^1 (X,
L^2(-x)) \bigcap \pp H^1 (X, L^{-2}(-x))$ when we consider their
images inside the moduli space $M_x$. For each $x \in X$, any
Hecke cycle in $M_X$ lying over $l \in \fK - \fK_0$ restricts on
$M_x$ to this kind of line pairs.

Next consider $L \in Pic^0(X)$ with $L^2 \cong \cO_X$ and let $l =
[L \oplus L] \in \fK_0$ be a nodal point. The fiber $\varphi^{-1}
(l)$ consists of two components $Q_l \cup R_l$: $Q_l$ is the space
of all conics which are contained in $\pp H^1 (X, \cO)$ and $R_l$
is the space of $\cO_{\pp^1}(-1)$-thickenings of lines in $\pp H^1
(X, \cO)$ which are contained in the thickening of $\pp H^1 (X,
\cO)_t$ (see \cite{NR} \S 3 and \S 4 for the details). The first
variety is isomorphic to a $\pp^5$-bundle over $Gr(\pp^2 ,
\pp^{g-1})=Gr(3,g)$ of planes in $\pp H^1 (X, \cO)$ while the
second variety is a $\pp^{g-2}$-bundle over the Grassmannian
$Gr(\pp^1, \pp^{g-1})=Gr(2,g)$ of lines in $\pp H^1 (X, \cO)$.

Finally we note that the fine moduli space $\bN$ of Hecke cycles
in $M_X$, has the following universal properties.

\begin{proposition}
(1) Suppose that there is a  flat family of closed subschemes of
$M_X$,
$$\xymatrix{
\mathcal{C}\ar@{^(->}[rr]\ar[dr] && M_X\times T\ar[dl]\\
&T}$$ parameterized by $T$ such that the fiber $\mathcal{C}_t$ is
a good Hecke cycle for generic $t \in T$. Then we have an induced
morphism $\tau : T \rightarrow \bN$ such that $\tau (t) =
[\mathcal{C}_t] \in \bN$.

(2) Suppose a holomorphic map $\tau: T \rightarrow \bN$ is given.
Suppose $T$ is an open subset of a nonsingular quasi-projective
variety $W$ on which a reductive group $G$ acts such that every
points in $W$ is stable and the smooth geometric quotient $W/G$
exists. Furthermore, assume that there is an open dense subset
$W'$ of $W$ such that whenever $t_1, t_2 \in T \cap W'$ are in the
same orbit, we have $\tau(t_1) = \tau(t_2)$. Then $\tau$ factors
through the image $\bar{T}$ of $T$ in the quotient $W/G$.
\end{proposition}
\begin{proof}
These are consequences of the universal property of Hilbert scheme
and GIT quotients.
\end{proof}


\section{Kirwan's desingularization}

In this section, we review the Kirwan's desingularization $\bK$.
Main reference is \cite{k5} and we also refer the reader to
\cite{kiem} for an explicit description of the desingularization
process for the case of genus 3 curves.

As we noted before,
$$
M_0 = M_0^s \sqcup (\fK - \fK_0) \sqcup \fK_0,
$$
where $\fK_0$ consists of the $2^{2g}$ nodal points in $\fK$.
Kirwan's desingularization $\bK$ is obtained as a result of
systematic blow-ups of $M_0$. Let $M_1$ be the blow-up of $M_0$
along the deepest strata $\fK_0$. By blowing up $M_1$ along the
proper transform of of the middle stratum $\fK$, we get Kirwan's
{\it partial} desingularization $M_2$. By taking one more blow-up
along the singular locus of $M_2$, we get the {\it full}
desingularization $\bK$.

The moduli space $M_0$ is constructed as the GIT quotient $\fR
\git G$, where $G = SL(p)$ and $\fR$ is a smooth quasi-projective
variety which is a subset of the space of holomorphic maps from
$X$ to the Grassmannian $Gr(2, p)$ of $2$-dimensional quotients of
$\cc^p$ where $p$ is a large even number.

Let $l \in \fK_0$ represent $L \oplus L^{-1}$ where $L^2 \cong
\cO_X$. There is a unique closed orbit in $\fR^{ss}$ lying over
$h$. By deformation theory, the normal space of this orbit is
$$
H^1(End_0(L \oplus L^{-1})) \cong H^1(\cO) \otimes sl(2)
$$
where the subscript $0$ denotes the trace-free part. By Luna's
slice theorem, there is a neighborhood of $l$ homeomorphic to
$(H^1(\cO) \otimes sl(2)) \git SL(2)$ since the stabilizer of $h$
is $SL(2)$ (\cite{k5}, (3.3)). More precisely, there is an
$SL(2)$-invariant locally closed subvariety $W$ in $\fR^{ss}$
containing $l$ and an $SL(2)$-invariant morphism $W \rightarrow
H^1(\cO) \otimes sl(2)$, \'{e}tale at $h$, such that we have the
following commutative diagram with all horizontal morphisms being
\'{e}tale.
\begin{equation}
\begin{CD}
G \times_{SL(2)}(H^1(\cO) \otimes sl(2)) @<<< G \times_{SL(2)}W
@>>> \fR^{ss} \\
  @VVV         @VVV             @VVV \\
(H^1(\cO \otimes sl(2)) \git SL(2)  @<<<  W \git SL(2)  @>>> M_0 \\
\end{CD}
\end{equation}
\vskip 10pt

Next, let $l \in \fK - \fK_0$ represent $L \oplus L^{-1}$ with
$L^2 \ncong \cO$. The normal space to the unique closed orbit over
$l$ is isomorphic to
$$
H^1(End_0(L \oplus L^{-1})) \cong H^1(\cO) \oplus H^1(L^2) \oplus
H^1(L^{-2}).
$$
Here the stabilizer $\cc^*$ acts with weights $0, 2, -2$
respectively on the components, and there is a neighborhood of $l$
homeomorphic to
$$
H^1(\cO) \bigoplus (H^1(L^2) \oplus H^1 (L^{-2}) \git \cc^*).
$$
Notice that $H^1(\cO)$ is the tangent space to $\fK$ and hence
$$
H^1(L^2) \oplus H^1 (L^{-2}) \git \cc^* \cong \cc^{2g-2} \git
\cc^*
$$
is the normal cone. The GIT quotient of the projectivization $\pp
\cc^{2g-2}$ by the induced $\cc^*$-action is $\pp^{g-2} \times
\pp^{g-2}$ and the normal cone $\cc^{2g-2} \git \cc^*$ is obtained
by collapsing the zero section of the line bundle $\cO_{\pp^{g-2}
\times \pp^{g-2}} (-1, -1)$.

Let $Z_{SL(2)}^{ss}$ (resp. $Z_{\cc^*}^{ss}$) be the set of
semiatable points in $\fR^{ss}$ fixed by $SL(2)$ (resp. $\cc^*$).
Let $\fR_1$ be the blow-up of $\fR^{ss}$ along the smooth
subvariety $GZ_{SL(2)}^{ss}$. Then by \cite{k2} Lemma 3.11, the
GIT quotient $\fR_1^{ss} \git G$ is the first blow-up $M_1$ of
$M_0$ along $G Z_{SL(2)}^{ss} \git G \cong \fK_0$. The
$\cc^*$-fixed point set in $\fR_1^{ss}$ is the proper transform
$\tilde{Z}_{\cc^*}^{ss}$ of $Z_{\cc^*}^{ss}$ and the quotient of
$G \tilde{Z}_{\cc^*}^{ss}$ by $G$ is the blow-up $\tilde{\fK}$ of
$\fK$ along $\fK_0$. Let $\fR_2$ be the blow-up of $\fR_1^{ss}$
along the smooth subvariety $G\tilde{Z}_{\cc^*}^{ss} = G
\times_{N^{\cc^*}} \tilde{Z}_{\cc^*}^{ss}$, where $N^{\cc^*}$ is
the normalizer of $\cc^*$. Then again by \cite{k2} Lemma 3.11, the
GIT quotient $\fR_2^{ss} \git G$ is the second blow-up $M_2$ of
$M_1$ along $G \tilde{Z}_{\cc^*}^{ss} \git G \cong \tilde{\fK}$.
This is Kirwan's partial desingularization of $M_0$ (see \S 3 of
\cite{k5}).

The points with stabilizer greater than the center $\{\pm 1 \}$ in
$\fR_2^{ss}$ is precisely the exceptional divisor of the second
blow-up and the proper transform $\tilde{\Delta}$ of the subset
$\Delta$ of the exceptional divisor of the first blow-up, which
corresponds, via Luna's slice theorem, to
$$
SL(2) \cdot \pp \{ \left(\begin{matrix} 0  & b \\ c & 0
\end{matrix}\right) | b,c \in H^1(\cO) \} \ \subset \ \pp
(H^1(\cO)\otimes sl(2)).
$$
Hence by blowing up $M_2$ along $\tilde{\Delta} \git SL(2)$, we
get a smooth variety $\bK$, Kirwan's desingularization.

We can now state the main result of this paper. Note that both
$\bN$ and $\bK$ contain $M_0^s$ as dense open subsets. Hence we
have a birational map $\rho':\bK\dashrightarrow \bN$.
\begin{theorem}\label{mainthm}
$\rho'$ extends to a morphism $\rho:\bK\to \bN$.\end{theorem}

In the subsequent two sections, we prove this theorem and in
section \ref{down} we show that $\rho$ is in fact a blow-up along
a smooth subvariety of $\bN$. Finally, in section \ref{comp8} we
compute the cohomology of $\bN$.


\section{Middle stratum}\label{midsec}
Let us first extend $\rho'$ to points over the middle stratum of
$M_0$. Let $l=[L\oplus L^{-1}]\in \mathfrak{K}-\fK_0$ be a
non-nodal point in the Kummer variety and let $W$ be the \'etale
slice of the unique closed orbit in $\mathfrak{R}^{ss}$ over $l$.
The deformation space of $L\oplus L^{-1}$ with determinant fixed
is
\begin{equation}\label{decNwt}
\cN=H^1(\End_0(L\oplus L^{-1}))=H^1(\cO)\oplus H^1(L^2)\oplus
H^1(L^{-2})\end{equation}
 where the subscript $0$ above denotes
the trace-free part. There is a versal deformation $\cF$ over
$\cN\times X$ and this gives us an analytic isomorphism of a
neighborhood $U$ of $0$ in $\cN$ with a neighborhood of $l$ in
$W$. The restriction of $\cF$ to $H^1(\cO)$ is a direct sum
$\cL\oplus \cL^{-1}$ where $\cL$ is the versal deformation of the
line bundle $L$. The group $\cc^*$ acts with weights $0,2,-2$
respectively on the three components of $\cN$ in \eqref{decNwt}.

Let $\pi:\tilde{\cN}\to \cN$ be the blow-up along $H^1(\cO)$ and
let $\tcN^{s}$ be the set of stable points in $\tcN$ with respect
to the obvious induced action of $\cc^*$. Let $D$ be the set of
stable points in the exceptional divisor of the blow-up; let
$\tcF$ be the pull-back of $\cF$ to $\tcN^s\times X$; let $\tcL$
be the pull-back of $\cL$ to $D$; let $\psi:\pp \tcF^*\to
\tcN^s\times X$ be the projectivization of the dual of $\tcF$.
Consider the composition
$$\xymatrix{
\pp \tcF^*\times X \ar[r]^(.4){\psi\times 1_X} & (\tcN^s\times
X)\times X\ar[r]^(.6){p_{13}} & \tcN^s\times X }$$ where $p_{13}$
denotes the projection onto the product of the first and the third
components. Let $\tcF'$ be the pull-back of $\tcF$ via the above
composition; let $q_X$ (resp. $q_N$) be the composition of $\psi$
with the projection onto $X$ (resp. $\tcN^s$); let $i:\pp
\tcF^*\to \pp \tcF^*\times X$ be the map $1_{\pp \tcF^*}\times
q_X$. Then there is a tautological homomorphism $\tcF'\to
i_*\cO_{\pp \tcF^*}(1)$. Let $\cE$ be its kernel. Then $\cE$ is a
family of rank 2 bundles on $X$ of degree $-1$ parameterized by
$\pp \tcF^*$ since for each $\theta\in \pp \tcF^*$,
\begin{equation}\label{defEmid} \cE|_{\{\theta\}\times
X}=\mathrm{ker}(\tcF|_{\{q_N(\theta)\}\times X}\twoheadrightarrow
\cO_{q_X(\theta)}).\end{equation}

The isomorphism $\tcF|_{D\times X}\cong \tcL\oplus \tcL^{-1}$
gives rise to two sections
$$s_1,s_2:D\times X\to \pp \tcF^*|_{D\times X}$$
by considering the obvious surjections $\tcF|_{D\times X}\to \tcL$
and $\tcF|_{D\times X}\to \tcL^{-1}$. Thus we have two disjoint
codimension 2 subvarieties $s_1(D\times X)$ and $s_2(D\times X)$
of $\pp \tcF^*$.

\begin{lemma}\label{lem6.1} $\cE|_{\{\theta\}\times X}$ is stable if and only if
$\theta\in \pp \tcF^*-s_1(D\times X)-s_2(D\times X)$.\end{lemma}
\begin{proof} If $q_N(\theta)\notin D$, $\cF|_{\{q_N(\theta)\}\times
X}$ is a stable bundle and hence $\cE|_{\{\theta\}\times X}$ is
stable since it is the result of an elementary modification at one
point (\cite{NR} Lemma 5.5). For $\theta \in
q_N^{-1}(D)-s_1(D\times X)-s_2(D\times X)$,
$\cE|_{\{\theta\}\times X}$ is the result of an elementary
modification $L_{\theta}\oplus L_{\theta}^{-1}\twoheadrightarrow
\cc$ for $L_{\theta}=\tcL|_{\{q_N(\theta)\}\times X}$ where
$L_{\theta}\to L_{\theta}\oplus L_{\theta}^{-1}\to \cc$ and
$L_{\theta}^{-1}\to L_{\theta}\oplus L_{\theta}^{-1}\to \cc$ are
both nonzero. It is an elementary exercise to show that the result
of this modification is a stable bundle, whose isomorphism class
is independent of the choice of the map $L_{\theta}\oplus
L_{\theta}^{-1}\twoheadrightarrow \cc$.

If $\theta\in s_1(D\times X)$, $\cE|_{\{\theta\}\times X}$ is
$L_{\theta}(-x)\oplus L_{\theta}^{-1}$ where $x=q_X(\theta)$.
Hence it is unstable. Similarly $\cE|_{\{\theta\}\times X}$ is
unstable for $\theta\in s_2(D\times X)$. \end{proof}

From the definition \eqref{defEmid}, we have $$\cE|_{s_1(D\times
X)\times X}\cong \tcL^{-1}\oplus \tcL(-\mathbf{q}_X)$$ where
$\tcL(-\mathbf{q}_X)$ is the kernel of $\tcL\to
\tcL|_{i(s_1(D\times X))}$. Similarly
$$\cE|_{s_2(D\times X)\times X}\cong \tcL\oplus
\tcL^{-1}(-\mathbf{q}_X).$$ In order to get a family of stable
bundles, we blow up $\pp \tcF^*$ along the locus of unstable
bundles $s_1(D\times X)\cup s_2(D\times X)$. Let $p:Z\to\pp
\tcF^*$ be this blow-up and $D', D''$ be the exceptional divisors
for $s_1$ and $s_2$ respectively. Let $\cE'$ denote the pull-back
of $\cE$ to $Z\times X$ and $\cL'$ and $\cL'(-\mathbf{q}_X)$
(resp. $\cL''$ and $\cL''(-\mathbf{q}_X)$) be the pull-backs of
$\tcL$ and $\tcL(-\mathbf{q}_X)$ to $D'\times X$ (resp. $D''\times
X$). Then we have $\cE'|_{D'\times X}\cong \cL'^{-1}\oplus
\cL'(-\mathbf{q}_X)$ and $\cE'|_{D''\times X}\cong \cL''\oplus
\cL''^{-1}(-\mathbf{q}_X)$. Now let
$$\bE=\mathrm{ker}\left[\cE' \to \cE'|_{(D' \cup D'') \times X}
\to \cL'(-\mathbf{q}_X)\oplus
\cL''^{-1}(-\mathbf{q}_X)\right].$$

\begin{lemma}\label{lem6.2}
$\bE$ is a family of stable vector bundles of degree $-1$ on $X$
parameterized by $Z$.\end{lemma}
\begin{proof}
Let $\theta$ be any point in $s_1(D\times X)$ and $x=q_X(\theta)$.
Let $C$ be a line in $\cN$ given by a map $\cc\to \cN$, $z\mapsto
(a, zb, zc)$ for $a\in H^1(\cO)$, $0\ne b\in H^1(L^2)$, $0\ne c\in
H^2(L^{-2})$. Note that any point in $D$ is represented by such a
line. We consider such a line for $q_N(\theta)$. By restricting to
a neighborhood $U$ of $0$ in $\cc$, we can find a finite open
covering $\{V_i\}$ of $X$ such that $\cF|_{U\times V_i}$ is
trivial and $x$ is contained only in $V_1$. Fix a trivialization
for each $i$. Then the transition matrix of
$F^z:=\cF|_{\{(a,zb,zc)\}\times X}$ from $V_i$ to $V_j$ is of the
form
\begin{equation}\label{transmid}
\left(\begin{matrix} \lambda_{ij} & zb_{ij}\\
zc_{ij}& \lambda_{ij}^{-1}\end{matrix}\right)\end{equation} where
$\{\lambda_{ij}|_{z=0}\}$ is the transition for
$L_a:=\cL|_{\{(a,0,0)\}\times X}$. Further, $b$ and $c$ are the
cohomology classes of the cocycles $\{\lambda_{ij}b_{ij}|_{z=0}\}$
and $\{\lambda_{ij}^{-1}c_{ij}|_{z=0}\}$ in $H^1(L_a^2)\cong
H^1(L^2)$ and $H^1(L_a^{-2})\cong H^1(L^{-2})$ respectively.

The normal space to $s_1(D\times X)$ at $\theta$ is a two
dimensional space $\{(z,t)\}$ where $(z,t)$ represents the bundle
$F^z$ and the surjection $F^z|_x\cong \cc^2\twoheadrightarrow \cc$
given by $(1,t)$. By definition, the bundle $\cE$ restricted to
$(z,t)\times X$ is the kernel of $F^z\to F^z|_x\to \cc$. Its
transition matrices can be described as follows. Let $\zeta$ be a
coordinate function on $V_1$ such that $\zeta(x)=0$. A section on
$V_1$ of the kernel is of the form $(\zeta f-tg, g)$ for some
holomorphic functions $f, g$. From the computation
\begin{equation}\label{eq5.4}
\left(\begin{matrix} f\\ g\end{matrix}\right)\leftrightarrow
\left(\begin{matrix} \zeta f-tg\\ g\end{matrix}\right) \mapsto
\left(\begin{matrix} \lambda_{1j} & zb_{1j}\\
zc_{1j}& \lambda_{1j}^{-1}\end{matrix}\right) \left(\begin{matrix}
\zeta f-tg\\ g\end{matrix}\right) = \left(\begin{matrix} \zeta
\lambda_{1j} & zb_{1j}-t\lambda_{1j}\\ \zeta z c_{1j} &
-ztc_{1j}+\lambda_{1j}^{-1}\end{matrix}\right)
\left(\begin{matrix} f\\ g\end{matrix}\right)\end{equation} the
transition matrix from $V_1$ to $V_j$ for $j\ne 1$ is
$$\left(\begin{matrix}
\zeta \lambda_{1j} & zb_{1j}-t\lambda_{1j}\\
\zeta z c_{1j} & -ztc_{1j}+\lambda_{1j}^{-1}\end{matrix}\right).$$
The transition from $V_j$ to $V_1$ for $j\ne 1$ is the inverse
matrix
$$\left(\begin{matrix}
\zeta^{-1} (\lambda_{j1}+ztc_{j1}) & \zeta^{-1}(zb_{j1}+t\lambda_{j1}^{-1})\\
z c_{j1} & \lambda_{j1}^{-1}\end{matrix}\right)$$ and the other
transition matrices are unchanged \eqref{transmid}.

Any point $\tilde{\theta}$ in $Z$ over $\theta$ is represented by
a line through $0$ in the $(z,t)$-plane. Suppose $t=\tau z$ for
some $\tau\in \cc$. When $z=0$ the transition matrices are
diagonal and the bundle is just $L_a(-x)\oplus L_a^{-1}$. To get
$\bE$, we modify $\cE$ by the surjection $\cE|_{(0,0)}\cong
L_a(-x)\oplus L_a^{-1}\to L_a(-x)$. A section over $V_1$ of $\bE$
restricted to the line $t=\tau z$ is of the form $(zf, g)$ for
some holomorphic functions $f,g$. From
\begin{equation}\label{eq5.5}\begin{array}{ll}
&\left(\begin{matrix} f\\ g\end{matrix}\right)\leftrightarrow
\left(\begin{matrix}z f\\ g\end{matrix}\right)\mapsto
\left(\begin{matrix}
\zeta \lambda_{1j} & zb_{1j}-\tau z\lambda_{1j}\\
\zeta z c_{1j} & -\tau
z^2c_{1j}+\lambda_{1j}^{-1}\end{matrix}\right)
\left(\begin{matrix}z f\\ g\end{matrix}\right)\\
&=\left(\begin{matrix} z(\zeta\lambda_{1j}f+(b_{1j}-\tau
\lambda_{1j})g)\\ \zeta z^2c_{1j}f+(-\tau z^2
c_{1j}+\lambda_{1j}^{-1})g\end{matrix}\right)\leftrightarrow
\left(\begin{matrix}
\zeta \lambda_{1j} & b_{1j}-\tau \lambda_{1j}\\
\zeta z^2c_{1j} & -\tau
z^2c_{1j}+\lambda_{1j}^{-1}\end{matrix}\right)
\left(\begin{matrix} f\\
g\end{matrix}\right)\end{array}
\end{equation}
we see that the transition matrix of $\bE|_{\tilde{\theta}\times
X}$ from $V_1$ to $V_j$ for $j\ne 1$ is
$$\left(\begin{matrix}
\zeta \lambda_{1j} & b_{1j}-\tau \lambda_{1j}\\
0 & \lambda_{1j}^{-1}\end{matrix}\right)$$ by plugging in $z=0$.
Similarly, the transition from $V_j$ to $V_1$ for $j\ne 1$ is
$$\left(\begin{matrix}
\zeta^{-1} \lambda_{j1} & \zeta^{-1}(b_{j1}+\tau \lambda_{j1}^{-1})\\
0 & \lambda_{j1}^{-1}\end{matrix}\right)$$ and the transition from
$V_i$ to $V_j$ for $i\ne 1$, $j\ne 1$ is
$$\left(\begin{matrix}
\lambda_{ij} & b_{ij}\\
0 & \lambda_{ij}^{-1}\end{matrix}\right).$$ This implies that
$\bE|_{\tilde{\theta}\times X}$ is an extension of $L_a^{-1}$ by
$L_a(-x)$. It is an elementary exercise to check that the
extension class in $H^1(L_a^2(-x))$ is given by
$$\mu^\tau_{ij}=\left\{\begin{array}{lll}
\lambda_{1j}(b_{1j}-\tau \lambda_{1j}) & \text{for} & i=1, j\ne 1\\
\zeta^{-1}\lambda_{i1}(b_{i1}+\tau \lambda_{i1}^{-1}) & \text{for}
& i\ne 1, j=1\\
\lambda_{ij}b_{ij} & \text{for} & i\ne 1, j\ne 1
\end{array}\right.$$ Note that
$$\mu^0_{ij}=\left\{\begin{array}{lll}
\lambda_{1j}b_{1j} & \text{for} & i=1, j\ne 1\\
\zeta^{-1}\lambda_{i1}b_{i1} & \text{for}
& i\ne 1, j=1\\
\lambda_{ij}b_{ij} & \text{for} & i\ne 1, j\ne
1\end{array}\right.$$ defines a class in $H^1(L_a^2(-x))$ which is
mapped to $b$ via the natural map $H^1(L_a^2(-x))\to H^1(L_a^2)$.
Similarly,
$$\mu^\infty_{ij}=\left\{\begin{array}{lll}
-\lambda_{1j}^2 & \text{for} & i=1, j\ne 1\\
\zeta^{-1} & \text{for}
& i\ne 1, j=1\\
0 & \text{for} & i\ne 1, j\ne 1\end{array}\right.$$ defines a
nonzero class in $H^1(L_a^2(-x))$ which generates the kernel of
$H^1(L_a^2(-x))\to H^1(L_a^2)$. Since $\mu_{ij}^\tau$ is a linear
combination of $\mu_{ij}^0$ and $\mu_{ij}^\infty$, the extension
classes for $\bE|_{p^{-1}(\theta)\times X}$ give us a projective
line in $\pp H^1(L_a^2(-x))$ which is also the projectivization of
the two dimensional subspace given by the preimage of $\cc b$.
Therefore $\bE|_{p^{-1}(\theta)\times X}$ is a family of stable
bundles. The same proof shows that $\bE|_{\tilde{\theta}\times X}$
is stable for $\tilde{\theta}\in D''$.
\end{proof}
\vskip 10 pt

\includegraphics{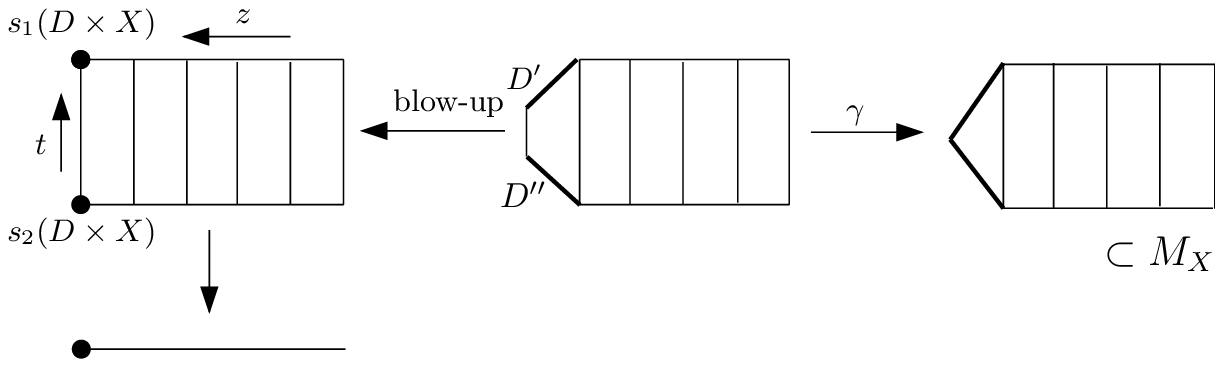}
\begin{center}
F{\small IGURE} 1
\end{center}

\vskip 10 pt

As a consequence of the above lemma, we get a morphism
$\gamma:Z\to M_X$ over $X$. By definition $\pp \tcF^*$ is a
projective line bundle over $\tcN^s\times X$ and hence flat over
$\tcN^s\times X$. For $\xi\in D$, $x\in X$, the fiber over
$(\xi,x) \in \pp \tcF^*$ in $Z$ is a chain of three rational
curves. As remarked in the proof of Lemma \ref{lem6.1} the
isomorphism class of the kernel of $L_a\oplus
L_a^{-1}\twoheadrightarrow \cc$ is independent of the surjection
if neither $L_a$ nor $L_a^{-1}$ is in the kernel. Hence $\gamma$
is constant on the middle component. The proof of Lemma
\ref{lem6.2} shows that the other two rational curves are embedded
by $\gamma$ into $\pp H^1(L_a^2(-x))$ and $\pp H^1(L_a^{-2}(-x))$
respectively as projective lines. By \cite{NR} Proposition 7.8,
this implies that the image of $p^{-1}(q_N^{-1}(\xi))$ by $\gamma$
is a limit Hecke cycle. On the other hand, for $\xi\in \tcN^s-D$,
the image of $p^{-1}(q_N^{-1}(\xi))$ by $\gamma$ is a good Hecke
cycle (Definition \ref{goodHecke}) and thus the Hilbert
polynomials of the fibers of the image by $\gamma$ of $Z$ over
$\tcN^s$ is constant. In particular, $\gamma(Z)$ is a flat family
of Hecke cycles in $M_X$ parameterized by $\tcN^s$. Therefore we
proved the following.
\begin{proposition}
 There is an analytic extension $\rho_l:\tcN^s\to \bN$ of
the obvious map $\rho'_l:\pi^{-1}(\cN^s)\to \bN$ which assigns
each stable bundle its associated good Hecke cycle where $\bN$ is
the moduli space of Hecke cycles in $M_X$.
\end{proposition}

Since two isomorphic stable bundles give us the same good Hecke
cycles, $\rho_l$ is invariant under the action of $\cc^*$ on the
open dense subset $\pi^{-1}(\cN^s)$ and hence $\rho_l$ is
$\cc^*$-invariant everywhere. So we get an analytic map
$$\overline{\rho}_l:\tcN^s\git \cc^*\to \bN .$$
Since a neighborhood of the vertex of the cone $\tcN^s\git \cc^*$
is analytically isomorphic to a neighborhood of $l\in
\mathfrak{K}-\fK_0\subset M_0$, we deduce that
$\rho:\bK\dashrightarrow \bN$ extends to the middle stratum
analytically.

\section{Deepest strata}\label{deepst}
In this section, we extend $\rho'$ to the points in $\bK$ over the
deepest strata $\fK_0=\zz_2^{2g}$. Since the exactly same argument
applies to every point in $\fK_0$, we consider only the points in
$\bK$ over $[\cO_X\oplus \cO_X]\in M_0$. The deformation space of
$\cO_X\oplus \cO_X$ with determinant fixed is
$$\cN=H^1(\cO_X)\otimes sl(2)$$ on which $SL(2)$ acts by conjugation
on $sl(2)$. There is a versal deformation $\cF$ over $\cN\times X$
which gives us an analytic isomorphism of a neighborhood of the
image $\overline{0}$ of $0$ in $\cN\git SL(2)$ with a neighborhood
of $[\cO_X\oplus \cO_X]$ in $M_0$.

Let $\Sigma$ be the subset of $\cN$ defined by
$$SL(2)\{H^1(\cO_X)\otimes \left(\begin{matrix}
1&0\\0&-1\end{matrix}\right)\}$$ which corresponds to the middle
stratum of $M_0$. Let $\pi_1:\cN_1\to \cN$ be the first blow-up in
the partial desingularization process, i.e. the blow-up at 0, and
let $\cD^{(1)}_1$ be the exceptional divisor. Let $\Delta$ be the
subset of $\cD^{(1)}_1$ defined as
$$SL(2)\pp \{\left(\begin{matrix} 0 & b\\
c&0\end{matrix}\right)\,|\, b,c\in H^1(\cO_X)\}$$ and let
$\tilde{\Sigma}$ be the proper transform of $\Sigma$ in $\cN_1$.
Then the singular locus of $\cN^{ss}_1\git SL(2)$ is the quotient
of $\Delta\cup \tilde{\Sigma}$ by $SL(2)$. It is an elementary
exercise to check that
\begin{equation}\label{eq4.0}\cD^{(1)}_1\cap
\tilde{\Sigma}=SL(2)\pp \{H^1(\cO_X)\otimes \left(\begin{matrix}
1&0\\0&-1\end{matrix}\right)\}=\Delta\cap
\tilde{\Sigma}.\end{equation} If we remove unstable points that
should be deleted after the desingularization process, $\Delta$ is
the locus in $\cD^{(1)}_1$ of $2\times 2$ matrices
$$\left(\begin{matrix} a & b\\
c&-a\end{matrix}\right)$$ with $\dim \mathrm{Span}\{a,b,c\}\le 2$
while $\Delta\cap \tilde{\Sigma}$ is the locus with $\dim
\mathrm{Span}\{a,b,c\}\le 1$.

 Let $\pi_2:\cN_2\to
\cN_1$ be the second blow-up, i.e. the blow-up along
$\tilde{\Sigma}$ and let $\cD^{(2)}_2$ be the exceptional divisor.
Let $\cD^{(1)}_2$ be the proper transform of $\cD^{(1)}_1$. The
singular locus of $\cN_2\git SL(2)$ is the quotient of the proper
transform $\tilde{\Delta}$ of $\Delta$.

Finally let $\pi_3:\tcN=\cN_3\to \cN_2$ denote the blow-up of
$\cN_2$ along $\tilde{\Delta}$ and let $\tcD^{(3)}=\cD^{(3)}_3$ be
the exceptional divisor while $\tcD^{(1)}=\cD^{(1)}_3$,
$\tcD^{(2)}=\cD^{(2)}_3$ are the proper transforms of
$\cD^{(1)}_2$ and $\cD^{(2)}_2$ respectively. Let $\pi:\tcN\to
\cN$ be the composition of the three blow-ups. Also let
$D^{(j)}_i$ be the quotient of $\cD^{(j)}_i$ in $\cN_i\git SL(2)$
for $1\le i\le 3$ and $1\le j\le i$.

\subsection{Modification over $\cN_1$}\label{subsecmodN1}
Let $\cN_1^{ss}$ be the set of semistable points in $\cN_1$. Let
$\cF_1$ be the pull-back of $\cF$ to $\cN^{ss}_1\times X$ and let
$\psi_1:\pp\cF_1^*\to \cN_1^{ss}\times X$ be the projectivization
of $\cF_1^*$. Consider the composition
$$\xymatrix{
\pp \cF_1^*\times X \ar[r]^(.4){\psi_1\times 1_X} &
(\cN_1^{ss}\times X)\times X\ar[r]^(.6){p_{13}} & \cN_1^{ss}\times
X }$$ where $p_{13}$ denotes the projection onto the product of
the first and the third components. Let $\cF_1'$ be the pull-back
of $\cF_1$ via the above composition; let $q_X$ (resp. $q_N$) be
the composition of $\psi_1$ with the projection onto $X$ (resp.
$\cN_1^{ss}$); let $i:\pp \cF_1^*\to \pp \cF_1^*\times X$ be the
map $1_{\pp \cF_1^*}\times q_X$. Then there is a tautological
homomorphism $\cF_1'\to i_*\cO_{\pp \cF_1^*}(1)$. Let $\cE_1$ be
its kernel. Then $\cE_1$ is a family of rank 2 bundles on $X$ of
degree $-1$ parameterized by $\pp \cF_1^*$. For $\theta_1\in
q_N^{-1}(\cDoo)$, $\cF_1'|_{\theta_1\times X}\cong \cO\oplus \cO$
and $\cE_1|_{\theta_1\times X}\cong \cO(-q_X(\theta_1))\oplus \cO$
which is unstable. We modify $\cE_1$ to get a family of stable
bundles on $\cDoo-\tsig$.

Since $\cDoo\subset \pi_1^{-1}(0)$, $\cF_1|_{\cDoo\times X}\cong
\cO\oplus \cO$ and hence $q_N^{-1}(\cDoo)=\pp
\cF^*_1|_{\cDoo\times X}=\pp^1\times \cDoo\times X.$ The
restriction $\cF'_1|_{q_N^{-1}(\cDoo)\times X}$ is thus $\cO\oplus
\cO$ and the tautological homomorphism $\cF_1'\to i_*\cO_{\pp
\cF_1^*}(1)$ restricted to $q_N^{-1}(\cDoo)\times X$ can be
factored as
$$\cF'_1|_{q_N^{-1}(\cDoo)\times X}\cong \cO\oplus \cO\to
\cO(1)\to \cO(1)|_{i(q_N^{-1}(\cDoo) )}$$ where $\cO(1)$ denotes
the pull-back of $\cO_{\pp^1}(1)$ by the projection
$q_N^{-1}(\cDoo)\times X\to \pp^1$. Let $\cO(-\mathbf{q}_X)$
denote the kernel of the above surjection $\cO(1)\to
\cO(1)|_{i(q_N^{-1}(\cDoo) )}$ over $q_N^{-1}(\cDoo)\times X$. By
definition, the composition $\cE_1|_{q_N^{-1}(\cDoo)\times X}\to
\cF'_1|_{q_N^{-1}(\cDoo)\times X}\to \cO(1)|_{i(q_N^{-1}(\cDoo)
)}$ is zero and thus we have a homomorphism
$$\cE_1\to \cE_1|_{q_N^{-1}(\cDoo)\times X}\to \cO(-\mathbf{q}_X).$$
Let $\bE_1$ be its kernel.

\begin{lemma}\label{lem7.1}
Let $\xi_1=\left[\begin{matrix} a&b\\c&-a\end{matrix}\right]\in
\cDoo=\pp \cN$ with $a,b,c\in H^1(\cO_X)$.\begin{enumerate}
\item Suppose $\dim \mathrm{Span}\{a,b,c\}=3$. Then
$\bE_1|_{q_N^{-1}(\xi_1)\times X}$  is a family of stable bundles
which gives us a morphism $\gamma_{\xi_1}:q_N^{-1}(\xi_1)\to M_X$
over $X$. Furthermore, the image of $\psi_1^{-1}(\xi_1,x)$ by
$\gamma_{\xi_1}$ for any $x\in X$ is a nonsingular conic in $\pp
H^1(\cO_X)\cong \pp H^1(\cO_X(-x))\subset M_x$.
\item For $\xi_1\in \Delta-\tsig$, $\bE_1|_{q_N^{-1}(\xi_1)\times
X}$ is a family of stable bundles and the map $\pp^1\cong
\psi_1^{-1}(\xi_1,x)\to M_x$ is a branched double covering onto a
projective line in $\pp H^1(\cO_X)\cong \pp H^1(\cO_X(-x))\subset
M_x$.
\end{enumerate}
\end{lemma}
\begin{proof}
We use the same method as in Lemma \ref{lem6.2}. Let $x\in X$ be
any. The line $\cc\to \cN$, $z\mapsto \left(\begin{matrix}
za&zb\\zc&-za\end{matrix}\right)$ represents $\xi_1$. By
restricting to a neighborhood $U$ of $0$ in $\cc$, we can find a
finite open covering $\{V_i\}$ of $X$ such that $\cF|_{U\times
V_i}$ is trivial and $x$ is contained only in $V_1$. Fix a
trivialization for each $i$. The transition matrix of
$F^z:=\cF|_{\{(za,zb,zc)\}\times X}$ from $V_i$ to $V_j$ is of the
form
\begin{equation}\label{deeptrans}
\left(\begin{matrix} 1+za_{ij} & zb_{ij}\\
zc_{ij}& 1-za_{ij} \end{matrix}\right)\end{equation} mod $z^2$.
Then $\{b_{ij}|_{z=0}\}$ and $\{c_{ij}|_{z=0}\}$ are cocycles
represented by $b,c\in H^1(\cO_X)$ respectively. The fiber $\pp^1$
over $(\xi_1,x)$ in $\pp \cF_1^*$ has two charts given by
$$(1,t):\cc^2\to \cc\qquad \qquad (s,1):\cc^2\to \cc.$$ Let us
consider the first chart $(1,t)$.

By definition, the bundle $\cE_1|_{(\xi_1,x,t)}$ is obtained as a
consequence of the elementary modification of $F^z|_{z=0}$ at $x$
by $(1,t):\cc^2\to \cc$. Let $E_1^z$ be the kernel of $F^z\to
F^z|_x\cong \cc^2\to \cc$ where the last map is $(1,t)$ and let
$\zeta$ be a coordinate function of $V_1$ with $\zeta(x)=0$. The
computation \eqref{eq5.4} tells us that the transition matrix of
$E_1^z$ from $V_1$ to $V_j$ for $j\ne 1$ is
$$A_{1j}=\left(\begin{matrix}
\zeta(1+za_{1j})& zb_{1j}-t(1+za_{1j})\\
\zeta zc_{1j} & -ztc_{1j}+1-za_{1j}\end{matrix}\right).$$ The
transition from $V_j$ to $V_1$ is the inverse matrix
$A_{j1}=A_{1j}^{-1}$ and the other transition matrices are
unchanged \eqref{deeptrans}.

Next $\bE_1|_{(\xi_1,x,t)}$ is the result of an elementary
modification at $z=0$ of the family $E_1=\{E^z_1\}\to \cc\times X$
parameterized by $\cc$. The transition matrix of $E_1^0$ from
$V_1$ to $V_j$ for
$j\ne 1$ is $$A_{1j}^0=\left(\begin{matrix} \zeta & -t\\
0&1\end{matrix}\right)$$ Consider the commutative diagram
\begin{equation}\label{dia1}
\xymatrix{ \cc^2\ar[d]_{(1,0)}\ar[r]^{A_{1j}^0} &\cc^2\ar[d]^{(1,t)}\\
\cc\ar[r]_{\zeta} &\cc
 }\end{equation}
The horizontal maps are the transitions from $V_1$ to $V_j$ for
$E_1^0$ and $\cO_X(-x)$ respectively. The transitions from $V_j$
to $V_1$ is the inverse matrices and the other transitions from
$V_i$ to $V_j$ ($i,j\ne 1$) are identity. The vertical maps, which
is $(1,0)$ for $V_1$ and $(1,t)$ for $V_i$, $i\ne 1$, give us the
surjection $E_1^0\twoheadrightarrow \cO_X(-x)$ and let
$\{\bE_1^z\}$ be the kernel of $$E_1\to
E_1|_{z=0}=E_1^0\twoheadrightarrow\cO_X(-x).$$ Then $\bE_1^0$ is
our $\bE_1|_{(\xi_1,x,t)}$. Let us find the transition matrices of
$\bE_1^0$. From \eqref{dia1}, a section of the kernel of $E_1\to
\cO_X(-x)$ on $V_1$ is of the form $(zf,g)$ for some holomorphic
functions $f$ and $g$. Also a section of the kernel on $V_j$ is of
the form $(zf-tg, g)$. Note that to recover $(f,g)$ from $(zf-tg,
g)$ we need to multiply $\left(\begin{matrix} z^{-1}&z^{-1}t\\ 0 &
1\end{matrix}\right)$. From the computation
\begin{equation}
\left(\begin{matrix} f \\ g\end{matrix}\right) \leftrightarrow
\left(\begin{matrix} zf \\ g\end{matrix}\right) \mapsto A_{1j}
\left(\begin{matrix} zf \\ g\end{matrix}\right)\leftrightarrow
\left(\begin{matrix} z^{-1}&z^{-1}t\\ 0 & 1\end{matrix}\right)
A_{1j}\left(\begin{matrix} z&0\\ 0 &
1\end{matrix}\right) \left(\begin{matrix} f \\
g\end{matrix}\right)\end{equation} we deduce that the transition
matrix of $\bE_1^z$ from $V_1$ to $V_j$ ($j\ne 1$) is
\begin{equation}\label{trantotal}
\left(\begin{matrix} \zeta(1+za_{1j}+ztc_{1j})& b_{1j}-2ta_{1j}-t^2c_{1j}\\
\zeta z^2 c_{1j} & 1-z a_{1j}-ztc_{1j}\end{matrix}\right)
\end{equation}
and thus the transition matrix $\bE_1^0$ from $V_1$ to $V_j$
($j\ne 1$) is
$$\left(\begin{matrix} \zeta& b_{1j}-2ta_{1j}-t^2c_{1j}\\ 0 &
1\end{matrix}\right)$$ after plugging in $z=0$. The transition
matrix from $V_j$ to $V_1$ is its inverse
$$\left(\begin{matrix} \zeta^{-1}&-\zeta^{-1}( b_{1j}-2ta_{1j}-t^2c_{1j})\\ 0 &
1\end{matrix}\right)$$ and the transition from $V_i$ to $V_j$
($i,j\ne 1$) is by a similar computation
$$\left(\begin{matrix} 1& b_{ij}-2ta_{ij}-t^2c_{ij}\\ 0 &
1\end{matrix}\right)$$ This implies that $\bE_1|_{(\xi_1,x,t)}$ is
an extension of $\cO_X$ by $\cO_X(-x)$. Via the isomorphism
$H^1(\cO_X(-x))\cong H^1(\cO_X)$, the extension class is given by
\begin{equation}\label{6.4.1}\mu_{ij}^t=b_{ij}-2ta_{ij}-t^2c_{ij}\end{equation}
 and thus it is
$b-2ta-t^2c$ in $H^1(\cO_X)$. If we use $(s,1)$ as our chart on
$\pp^1$, we get the extension class $s^2b-2sa-c$ similarly.
Therefore, $\bE_1|_{\psi_1^{-1}(\xi_1,x)}$ gives us the locus
$\{s^2b-2sta-t^2c\,|\, [s,t]\in \pp^1\}$ in $\pp H^1(\cO_X)\cong
\pp H^1(\cO_X(-x))\hookrightarrow M_x$. If $a,b,c$ are
independent, the locus is a nonsingular conic.

The points in $\Delta$ are of the form $$\left[\begin{matrix}
0&b\\ c&0\end{matrix}\right]$$ after conjugation. In this case,
the above locus is a line in $\pp H^1(\cO_X)$ and the map $\pp^1
\cong \psi_1^{-1}(\xi_1,x)\to \pp H^1(\cO_X)$ is a branched double
covering.
\end{proof}

Note that if $\dim \mathrm{Span} \{ a,b,c\}=2$, the matrix
$\left[\begin{matrix} a&b\\ c&-a\end{matrix}\right]$ is conjugate
to matrices of the form $$\left[\begin{matrix} *&0\\
*&*\end{matrix}\right]\qquad\text{or}\qquad \left[\begin{matrix} 0&*\\
*&0\end{matrix}\right]$$ The first case becomes unstable in
$\cN_2$ and thus should be removed after all. The second matrix
lies in $\Delta-\tsig$. Hence the above lemma says $\bE_1$ is a
family of stable bundles when $\dim \mathrm{Span} \{ a,b,c\}\ge
2$. In the next subsection, we deal with the case when $\dim
\mathrm{Span} \{ a,b,c\}=1$.

\subsection{Modification over $\cN_2$}\label{subsN2}

Let $\cF_2$ be the pull-back of $\cF_1$ by $\pi_2\times 1:
\cN_2^s\times X\to \cN_1^{ss}\times X$ where $\cN_2^s=\cN_2^{ss}$
is the set of stable points in $\cN_2$. Let $\psi_2:\pp\cF_2^*\to
\cN_2^s\times X$ be the projectivization of $\cF_2^*$. By abuse of
notation, let $q_N$ (resp. $q_X$) denote the composition of
$\psi_2$ with the projection onto $\cN_2^s$ (resp. $X$). Then $\pp
\cF_2^*$ is the pull-back of $\pp \cF_1^*$. Let $\pp \cF_2^*\to
\pp \cF_1^*$ be the obvious map and let $\cE_2$ be the pull-back
of $\bE_1$ to $\pp \cF_2^*\times X$.

\begin{lemma} The locus of unstable bundles $S=\{\theta \in \pp
\cF_2^*\,|\,\cE_2|_{\theta\times X} \text{ is unstable}\}$ is a
smooth subvariety of codimension 2. Furthermore, $\cE_2|_{S\times
X}\cong \cL\oplus \cM$ where $\cL$ (resp. $\cM$) is a family of
line bundles of degree $0$ (resp. degree $-1$).
\end{lemma}
\begin{proof}
The modification of a semistable rank 2 bundle $F$ with
$\mathrm{det} F=\cO$ on $X$ by $F\to F|_x\cong
\cc^2\twoheadrightarrow\cc$ is unstable if and only if $F$ is an
extension $0\to L\to F\to L^{-1}\to 0$ for a line bundle $L$ of
degree 0 and the surjection $\cc^2\twoheadrightarrow\cc$ is
$F|_x\to L^{-1}|_x$. For $\xi\in \cN_2^s$, $\cF_2|_{\xi\times X}$
is a polystable bundle (because non-polystable bundles become
unstable in $\cN_2$) and the locus of strictly polystable bundles
in $\cN_1^{ss}$ is $\tsig\cup \cDoo$. Hence $S$ lies over
$\cDto\cup \cDtt $. But by Lemma \ref{lem7.1},
$\bE_1|_{q_N^{-1}(\cDoo-\tsig)}$ is a family of stable bundles.
Hence in fact $S$ lies over $\cDtt=\pi_2^{-1}(\tsig) $.

The proof of Lemma \ref{lem7.1} says for
$\xi_1=\left[\begin{matrix} a&0\\ 0&-a\end{matrix}\right]\in
\tsig\cap \cDoo$ and $x\in X$, $\bE_1|_{\psi_1^{-1}(\xi_1,x)}$ is
a family of extensions of $\cO_X$ by $\cO_X(-x)$ which splits at
exactly two points $(1,0)$ and $(0,1)$. Note that any point in
$\tsig\cap \cDoo$ is conjugate to $\left[\begin{matrix} a&0\\
0&-a\end{matrix}\right]$ for some $a\in H^1(\cO_X)$.

For $\xi_1\in \tsig$, $\cF_1|_{\xi_1\times X}$ is a direct sum of
line bundles $L\oplus L^{-1}$ for some line bundle $L$ of degree 0
and the locus of unstable bundles of $\bE_1$ in
$\psi_1^{-1}(\xi_1,x)=\pp \cF_1^*|_{(\xi_1,x)}\cong \pp^1$ for any
$x\in X$ is the two projections $L\oplus L^{-1}\twoheadrightarrow
L$ and $L\oplus L^{-1}\twoheadrightarrow L^{-1}$.

Let $J=\{\left(\begin{matrix} a&0\\
0&-a\end{matrix}\right)\,|\, a\in H^1(\cO_X)\}\subset \cN$. The
restriction of $\cF\to \cN\times X$ to $J\times X$ is $\cL\oplus
\cL^{-1}$ where $\cL$ is the versal deformation of the line bundle
$\cO_X$ over $H^1(\cO_X)\times X$ via the isomorphism $J\cong
H^1(\cO_X)$. Let $\tilde J$ be the blow-up of $J$ at 0 and $T\cong
\cc^*$ be the diagonal torus in $SL(2)$. Then the set of $T$-fixed
points in $\cN_1^{ss}$ is $\tilde J$ and $\tsig\cong SL(2)\times
_{N^T}\tilde J$ where $N^T$ is the normalizer of $T$ in $SL(2)$.
(cf.\cite{k5})

Consider the quotient map
$$SL(2)\times \tilde{J} \to SL(2)\times _{N^T}\tilde{J}\cong
\tsig.$$ The pull-back $\cF^{\dagger}$ of $\cF_1|_{\tsig\times X}$
using this map is isomorphic to the pull-back $\cF^{\sharp}$ of
$\cL\oplus \cL^{-1}$ by the projection $SL(2)\times \tilde{J}\to
\tilde{J}$. In fact, the isomorphism $\cF^\sharp\to \cF^\dagger$
over $(g,j)\in SL(2)\times \tilde{J}$ is given by
$$L_j\oplus L_j^{-1}\mapsto g(L_j\oplus L_j^{-1})g^{-1}.$$
Hence the two projections $\cL\oplus
\cL^{-1}\twoheadrightarrow\cL$ and $\cL\oplus
\cL^{-1}\twoheadrightarrow\cL^{-1}$ in $\pp (\cF^\sharp)^*$ give
us two sections of $\pp (\cF^\dagger)^*$. Note that if $g\in N^T$,
the union of the two sections is mapped to itself by conjugation
by $g$. Since the action of $N^T$ on $SL(2)\times \tilde{J}$ is
free, the union of the two sections descends to a smooth
subvariety of $\pp \cF_1^*|_{\tsig}$. Hence the locus of unstable
bundles of $\bE_1$ in $\pp\cF_1^*$ is a codimension 1 smooth
subvariety of $\pp \cF_1^*|_{\tsig}$. This implies that $S$ is a
codimension 2 subvariety of $\pp \cF_2^*$ lying over $\cDtt$.

For the second statement, let $\cE^\dagger$ be the kernel of the
tautological map from the pull-back of $\cF^\dagger$ to $\pp
(\cF^\dagger)^*\times X$ onto $i_*\cO_{\pp (\cF^\dagger)^*}(1)$
where $i:\pp (\cF^\dagger)^*\to \pp (\cF^\dagger)^*\times X$ is
$1_{\pp (\cF^\dagger)^*}\times q_X$, exactly as in the
construction of $\cE_1$ in subsection \ref{subsecmodN1}. Then it
is obvious from the isomorphism $\cF^\sharp\cong \cF^\dagger$ that
$\cE^\dagger$ restricted to the two sections is a direct sum of
line bundles of degree $0$ and $-1$ respectively. The action of
$N^T/T$ interchanges $L$ and $L^{-1}$. It also interchanges the
surjections $L\oplus L^{-1}\twoheadrightarrow L$ and $L\oplus
L^{-1}\twoheadrightarrow L^{-1}$. This implies that the line
bundles descend to $\tsig$ and hence we have the desired
decomposition of $\cE_2|_{S\times X}$.
\end{proof}

To remove unstable bundles from the family $\cE_2$ we proceed as
in section \ref{midsec}. Let $Z\to \pp\cF_2^*$ be the blow-up of
$\pp\cF_2^*$ along $S$; let $\cD$ be the exceptional divisor; let
$\cL'$ (resp. $\cM'$) be the pull-back of $\cL$ (resp. $\cM$) to
$\cD$; let $\cE'_2$ be the pull-back of $\cE_2$ to $Z\times X$;
let $\bE_2$ be the kernel of $$\cE_2'\to \cE_2'|_\cD\cong
\cL'\oplus \cM'\twoheadrightarrow \cM'.$$

\begin{lemma} $\bE_2$ is a family of rank 2 stable bundles of
degree $-1$.\end{lemma}
\begin{proof}
For the points over $\tsig-\Delta$, the proof is identical to
Lemma \ref{lem6.2}. For the points over $\tsig\cap \Delta$, we may
assume it lies over $\left[\begin{matrix} a& 0\\
0&-a\end{matrix}\right]$ for $a\in H^1(\cO_X)$ after conjugation.
The proof is then identical to that of Lemma \ref{lem6.2} if we
put $\lambda_{ij}=1$. The details are repetition of the same
computation and so we omit.
\end{proof}

Consequently, we have a morphism $$\gamma:Z\to M_X$$ over $X$. For
$\xi\in \cDtt$ and $x\in X$, the fiber over $(\xi,x)$ in $Z$ is a
chain of 3 rational curves. Since $b=c=0$ in \eqref{6.4.1}, the
extension class is a constant multiple of $a$, and therefore
$\gamma$ is constant on the middle component. As in section
\ref{midsec}, each of the other two rational curves is embedded
into $M_X$ by $\gamma$. When $\xi$ is not in the proper transform
$\tilde{\Delta}$ of $\Delta$ in $\cN_2$, the images of two curves
by $\gamma$ intersect transversely at one point and the image of
the fiber $(\xi,x)$ in $Z$ by $\gamma$ is a limit Hecke cycle as
in the middle stratum case. Therefore, we have a family of Hecke
cycles in $M_X$ parameterized by $\cN_2^s-\tilde{\Delta}$.

If $\xi\in \tilde{\Delta}$, the images of the two rational curves
by $\gamma$ coincide. To get Hecke cycles over $\tilde{\Delta}$ we
need to lift the family to $\tcN=\cN_3$.

\subsection{Hecke cycles over $\tdel$}\label{subsN3}
Recall that $\pi_3:\tcN=\cN_3\to \cN_2$ is the blow-up of $\cN_2$
along $\tdel$, $\tcD^{(3)}=\cD^{(3)}_3=\pi_3^{-1}(\tdel)$ and
$\tcN^s$ is the set of stable points in $\cN$. Let $\tilde{Z}$ be
the pull-back of $Z$ by $\pi_3\times 1_X$ so that we have the
diagram
$$\xymatrix{
\tilde{Z} \ar[r]^{\alpha}\ar[d]_{\tilde{\psi}} & Z\ar[d]^{\psi_2}\\
\tcN^s\times X \ar[r]^{\pi_3\times 1_X} & \cN_2^s\times X }$$ Let
$\tilde{\gamma}:\tilde{Z}\to Z\to M_X$ be the composition of
$\gamma$ with $\alpha$ and consider the diagram
$$\xymatrix{
\tilde{Z}\ar[rr]^{\tilde{\gamma}\times q}\ar[dr]_{q} &&M_X\times
\tcN^s\ar[dl]^{p_2}\\
&\tcN^s}$$ where $q:\tilde{Z}\to \tcN^s\times X\to \tcN^s$ is the
composition of $\tilde{\psi}$ with the projection onto $\tcN^s$
and $p_2$ is the projection onto the second component. Let
$\Gamma$ be the image of $\tilde{Z}$ by $\tilde{\gamma}\times q$
and $\phi$ be the restriction of $p_2$ to $\Gamma$. Then $\Gamma$
is a family of subschemes of $M_X$ parameterized by $\tcN^s$.

\begin{lemma} $\Gamma$ is a family of Hecke cycles.\end{lemma}
\begin{proof}
We have to show that the fiber $\Gamma_\xi:=\phi^{-1}(\xi)$ for
$\xi\in \tcD^{(3)}$ is a limit Hecke cycle. Every point in
$\tcD^{(3)}$ represents a normal direction of $\tdel$ in
$\cN_2^s$. After conjugation, we may assume $\xi_2:=\pi_3(\xi)\in
\tdel$ is of the form $$\left[\begin{matrix} 0&b\\
c&0\end{matrix}\right]$$ for some nonzero $b,c\in H^1(\cO_X)$. If
we restrict $\bE_2$ to the direction normal to $\cDto$ at $\xi_2$,
the transition matrix from $V_1$ to $V_j$ is given by
\eqref{trantotal} with $a_{1j}=0$ i.e.
$$\left(\begin{matrix}
\zeta(1+tzc_{1j}) & b_{1j}-t^2c_{1j}\\
\zeta z^2c_{1j} & 1-tzc_{1j}\end{matrix}\right)$$ mod $z^2$ and
the other transition matrices are given similarly. If $\xi$
represents a direction tangent to $\cDto$ at $\xi_2$, the
transition matrix of $\bE_2$ from $V_1$ to $V_j$ is
$$\left(\begin{matrix}
\zeta & b_{1j}-t^2c_{1j}-2tza_{1j}\\
0 & 1\end{matrix}\right)$$ by replacing $a_{1j}$ by $za_{1j}$ in
\eqref{trantotal} and the other transition matrices are given
similarly. In general, the normal direction represented by $\xi$
is a combination of the above two cases. Hence the transition from
$V_1$ to $V_j$ is
$$\left(\begin{matrix}
\zeta(1+tzc_{1j}) & b_{1j}-t^2c_{1j}-2tza_{1j}\\
\zeta z^2c_{1j} & 1-tzc_{1j}\end{matrix}\right)$$ and the other
transition matrices are given similarly. Thus the first order
variation in $z$ for the transition from $V_1$ to $V_j$ is
\begin{equation}\label{eq7.1a}
t \left(\begin{matrix} \zeta c_{1j}& -2a_{1j}\\
0 & -c_{1j}\end{matrix}\right)\end{equation} and those for the
other transitions are given similarly.

The image of $\psi_2^{-1}(\xi_2,x)$ for $x\in X$ is a projective
line $\pp^1$ in $\pp H^1(\cO_X(-x))$ from Lemma \ref{lem7.1} and
$t$ is a section of $\cO_{\pp^1}(1)$. Furthermore, the matrix in
\eqref{eq7.1a} represents a tangent vector in the moduli space of
``triangular bundles'' $\pp D$ over the Jacobian $Jac_0$ for $X$
of degree $0$, whose fiber over $L\in Jac_0$ is $\pp
H^1(L^2(-x))$. See \cite{NR} \S6.

Now observe that $\gamma$ is invariant under the $\zz_2$-action
given by $z\to -z$. In fact, the stabilizer of $\xi_2\in \tdel$ in
$SL(2)$ is $\zz_2\times \zz_2$. The first factor $\zz_2$ is the
center of $SL(2)$ and acts trivially everywhere. But the second
factor $\zz_2$ acts as $-1$ on the normal directions. Hence the
scheme theoretic fiber of $\Gamma$ over $\xi$ is the projective
line thickened by \eqref{eq7.1a}. This is more precisely the
thickening of $\pp^1$ by $\cO_{\pp^1}(-1)$ (since $t$ is a section
of $\cO_{\pp^1}(1)$) inside $\pp D$. By \cite{NR} Proposition 7.8,
the fiber $\Gamma_\xi$ is a limit Hecke cycle.
\end{proof}

By the above lemma, we have a map $\rho_0:\tcN^s\to \bN$. Since
$\rho_0$ is $SL(2)$-invariant on the dense open subset
$\pi^{-1}(\cN^s)$, it is invariant everywhere. Therefore, we have
a continuous map $$\overline{\rho}_0: \tcN^s/SL(2)\to \bN$$ which
implies that $\rho'$ extends to everywhere in $\bK$.


\section{Blowing down Kirwan's desingularization}\label{down}

Based on O'Grady's work \cite{ogrady}, it is shown in \cite{KL}
that $\bK$ can be blown down twice
\begin{equation}\label{eqf}\xymatrix{
f:&\bK\ar[r]^{f_{\sigma}}&\bK_{\sigma}\ar[r]^{f_{\epsilon}}&\bK_{\e}}.\end{equation}
Furthermore, they show in \cite{KL} that $\bK_\e$ is isomorphic to
Seshadri's desingularization of $M_0$ defined in \cite{Se1}. In
this section, we show that the moduli of Hecke cycles $\bN$ is in
fact the intermediate variety $\bK_\sigma$ which was conjectured
in \cite{KL}.

Let $\cA$ (resp. $\cB$) be the tautological rank 2 (resp. rank 3)
bundle over the Grassmannian $Gr(2,g)$ (resp. $Gr(3,g)$). Let
$W=sl(2)^{\vee}$ be the dual vector space of $sl(2)$. Fix $B\in
Gr(3,g)$. Then the variety of complete conics $\mathbf{CC} (B)$ is
the blow-up
$$\xymatrix{
\pp (S^2B)&\mathbf{CC}(B)\ar[l]_{\Phi_B}\ar[r]^{\Phi^{\vee}_B}&\pp
(S^2 B^{\vee})}$$ of both of the spaces of conics in $\pp B$ and
$\pp B^{\vee}$ along the locus of rank 1 conics. We recall the
following from \cite{KL} section 5.
\begin{proposition}\label{prop5.1}
\begin{enumerate}
\item $\tD^{(1)}$ is the variety of complete conics
$\mathbf{CC}(\cB)$ over $Gr(3,g)$. In other words, $\tD^{(1)}$ is
the blow-up of the projective bundle $\pp (S^2\cB)$ along the
locus of rank 1 conics.
\item There is an integer $l$ such that $$\tD^{(3)}\cong
\pp (S^2\cA)\times_{Gr(2,g)} \pp (\cc^g/\cA\oplus \cO(l)).$$ Hence
$\tD^{(3)}$ is a $\pp^2\times \pp^{g-2}$ bundle over $Gr(2,g)$.
\item The intersection $\tD^{(1)}\cap \tD^{(3)}$ is isomorphic to
the fibred product
$$\pp(S^2\cA)\times \pp (\cc^g/\cA)$$ over $Gr(2,g)$. As a
subvariety of $\tD^{(1)}$, $\tD^{(1)}\cap \tD^{(3)}$ is the
exceptional divisor of the blow-up $\mathbf{CC}(\cB)\to
{\pp}(S^2\cB^{\vee})$.
\item The intersection $\tD^{(1)}\cap\tD^{(2)}\cap \tD^{(3)}$ is isomorphic to
$$\pp(S^2\cA)_1\times \pp (\cc^g/\cA)$$ over $Gr(2,g)$ where $\pp(S^2\cA)_1$
denotes the locus of rank 1 quadratic forms.
\item The intersection $\tD^{(1)}\cap \tD^{(2)}$ is the
exceptional divisor of the blow-up $\mathbf{CC}(\cB)\to
{\pp}(S^2\cB)$.
\end{enumerate}\end{proposition}

Let $\s$ be the class of lines in the fiber of $\Phi_B^{\vee}$.
Then $\s$ gives us an extremal ray with respect to the canonical
bundle of $\bK$ and thus we can contract the ray. This turns out
to be the contraction of the $\pp(S^2\cA)$-direction of
$\tD^{(3)}$ and the contraction is a blow-down map $f_\s$. See
section 5 of \cite{KL} for details.

\begin{proposition}\label{iso7}
$\rho:\bK\to\bN$ factors through $\bK_\s$ and we have an
isomorphism $\bK_\s\cong \bN$.
\end{proposition}
\begin{proof}
By Riemann's extension theorem \cite{Kran}, it suffices to show
that $\rho$ is constant on the fibers of $f_{\s}$. From
Proposition \ref{prop5.1}, we know $f_{\s}$ is the result of
contracting the fibers $\pp^2$  of
$$\tD^{(3)}=\pp (S^2\cA)\times \pp (\cc^g/\cA\oplus \cO(l))\to
\pp (\cc^g/\cA\oplus \cO(l))$$ which amounts to forgetting the
choice of $b, c$ in the 2-dimensional subspace of $H^1(\cO)$
spanned by $b,c$. From our description of the transition matrices
of $\bE_2$ in subsection \ref{subsN2} and the thickening in
subsection \ref{subsN3}, it is easy to see that the Hecke cycles
on $\tcD^{(3)}$ depends on the two dimensional subspace spanned by
$\{b,c\}$ in $H^1(\cO_X)$ but not on the choices of $b,c$ in the
subspace. Hence $\rho$ factors through $\bK_\s$.

Now $\rho$ is an isomorphism over the stable part $M_0^s$ of
$M_0$. Further, the divisor $\tD^{(1)}$ is mapped to the divisor
$Q_k$ in \cite{NR} section 7.7 and the divisor $\tD^{(2)}$ is
mapped to the divisor given by $R_k$ in \cite{NR} 7.7. The
complements of (the images of) these sets in $\bK_\s$ and $\bN$
are of codimension $\ge 2$. Now by Zariski's main theorem, we
conclude that the induced map from $\bK_\s$ to $\bN$ is an
isomorphism.
\end{proof}
\begin{remark}
M.S. Narasimhan and S. Ramanan conjectured that the
desingularization $\bN$ can be blown down along certain projective
fibrations to obtain another nonsingular model of $M_0$(\cite{NR},
page 292) and this was proved by N. Nitsure(\cite{Ni}, Proposition
4.A.1 and 4.A.2). Our results, combined with \cite{KL}, show that
this blown-down process corresponds to the morphism
$$
f_\epsilon : \ \bK_\sigma (\cong \bN) \longrightarrow \bK_\epsilon
(\cong \bS).
$$
See \cite{KL} \S 5 for the structure of the morphism $f_\epsilon$.
\end{remark}

\section{Cohomology computation}\label{comp8}
In this section we compute the cohomology of the moduli of Hecke
cycles. For a variety $T$, let
$$P(T)=\sum_{k=0}^{\infty} t^k\dim H^k(T)$$ be the Poincar\'e
series of $T$. In \cite{k2}, Kirwan described an algorithm for the
Poincar\'e series of a partial desingularization of a good
quotient of a smooth projective variety and in \cite{k5} the
algorithm was applied to the moduli space without fixing the
determinant. For $P(M_2)$ we use Kirwan's algorithm in \cite{k2}.

Recall that $M_0=\mathfrak{R}^{ss}\git G$ where $G=SL(p)$ and
$\mathfrak{R}$ is a subset of the space of holomorphic maps from
$X$ to $Gr(2,p)$ for any sufficiently large even integer $p$
(\cite{k5} section 2). By \cite{AB} \S11 and \cite{KirL} \S13.1,
it is well-known that the equivariant Poincar\'e series
$P^G(\mathfrak{R}^{ss})=\sum_{k\ge 0}t^k \dim
H^k_G(\mathfrak{R}^{ss})$ is
$$\frac{(1+t^3)^{2g}-t^{2g+2}(1+t)^{2g}}{(1-t^2)(1-t^4)}+O(t^k)$$
where $k$ tends to infinity with $p$. Fix $p$ large enough so that
$k>6g-6$. In order to get $\mathfrak{R}^{ss}_1$ we blow up
$\mathfrak{R}^{ss}$ along $GZ^{ss}_{SL(2)}$ and delete the
unstable strata. So we get
$$P^G(\mathfrak{R}_1^{ss})=P^G(\mathfrak{R}^{ss})+2^{2g}\big(\frac{t^2+t^4+\cdots+t^{6g-2}}{1-t^4}
-\frac{t^{4g-2}(1+t^2+\cdots+t^{2g-2})}{1-t^2}\big).$$ Now
$\mathfrak{R}_2^{ss}$ is obtained by blowing up
$\mathfrak{R}_1^{ss}$ along $G\tilde{Z}^{ss}_{\cc^*}$ and deleting
the unstable strata. Thus we have
\begin{equation}\begin{array}{ll}
P^G(\mathfrak{R}_2^{ss})=P^G(\mathfrak{R}_1^{ss})
&+(t^2+t^4+\cdots+t^{4g-6})\big(\frac12
\frac{(1+t)^{2g}}{1-t^2}+\frac12
\frac{(1-t)^{2g}}{1+t^2} +2^{2g}\frac{t^2+\cdots+t^{2g-2}}{1-t^4}\big)\\
&-\frac{t^{2g-2}(1+t^2+\cdots+t^{2g-4})}{1-t^2}\big((1+t)^{2g}+2^{2g}(t^2+t^4+\cdots+t^{2g-2})\big).
\end{array}
\end{equation}
Because the stabilizers of the $G$ action on $\mathfrak{R}^{ss}_2$
are all finite, we have $$H^*_G(\mathfrak{R}_2^{ss})\cong
H^*(\mathfrak{R}^{ss}_2/G)=H^*(M_2)$$ and hence we deduce that
\begin{equation}\label{N2Pr}
\begin{array}{ll}
P(M_2)&=\frac{(1+t^3)^{2g}-t^{2g+2}(1+t)^{2g}}{(1-t^2)(1-t^4)}\\
&+2^{2g}\big(\frac{t^2+t^4+\cdots+t^{6g-2}}{1-t^4}
-\frac{t^{4g-2}(1+t^2+\cdots+t^{2g-2})}{1-t^2}\big)\\
&+(t^2+t^4+\cdots+t^{4g-6})\big(\frac12
\frac{(1+t)^{2g}}{1-t^2}+\frac12
\frac{(1-t)^{2g}}{1+t^2} +2^{2g}\frac{t^2+\cdots+t^{2g-2}}{1-t^4}\big)\\
&-\frac{t^{2g-2}(1+t^2+\cdots+t^{2g-4})}{1-t^2}\big((1+t)^{2g}+2^{2g}(t^2+t^4+\cdots+t^{2g-2})\big).
\end{array}
\end{equation}
Kirwan's desingularization is the blow-up of $M_2$ along
$\tilde{\Delta}\git SL(2)$ which is isomorphic to the $2^{2g}$
copies of  $\pp(S^2\cA)$ over $Gr(2,g)$. Hence,
$$P(\bK)=P(M_2)+2^{2g}(1+t^2+t^4)P(Gr(2,g))(t^2+t^4+\cdots + t^{2g-4})$$
by \cite{GH} p. 605.\footnote{The formula in \cite{GH} is stated
for smooth manifolds. But the same Mayer-Vietoris argument gives
us the same formula in our case (of orbifold $M_2$ blown up along
a smooth subvariety). The only thing to be checked is that the
pull-back homomorphism $H^*(M_2)\to H^*(\bK)$ is injective but
this is easy.}

On the other hand, $\bK$ is the blow-up of $\bK_{\s}$ along a
$\pp^{g-2}$-bundle over $Gr(2,g)$. Hence,
$$\begin{array}{ll}
P(\bN)=P(\bK_{\s})&=P(\bK)-2^{2g}(1+t^2+\cdots+t^{2g-4})P(Gr(2,g))(t^2+t^4)\\
&=P(M_2)+2^{2g}P(Gr(2,g))\frac{t^6-t^{2g-2}}{1-t^2}.\end{array}$$
By Schubert calculus \cite{GH}, we have
$$P(Gr(2,g))=\frac{(1-t^{2g})(1-t^{2g-2})}{(1-t^2)(1-t^4)}$$
and hence we proved the following.

\begin{proposition} The Poincar\'e polynomial of
$\bN$ is
\begin{eqnarray*}  P(\bN)
&=& \frac{(1+t^3)^{2g}}{(1-t^2)(1-t^4)}-
\frac{t^{2g-2}(1+t^2+t^4)(1+t)^{2g}}{(1-t^2)(1-t^4)}\\ && +
\frac{t^2}{2(1-t^2)}\Bigl[\frac{(1+t^{4g-6})(1+t)^{2g}}{1-t^2}
+\frac{(1-t^{4g-6})(1-t)^{2g}}{1+t^2}\Bigr] \\ && +
2^{2g}\Bigl[\frac{t^2(1-t^{6g-6})(1+t^4)}{(1-t^2)^2(1-t^4)}-
\frac{t^{2g-2}(1-t^{6})(1-t^{2g})(1+ t^{4})}{(1-t^2)^3(1-t^4)}
\Bigr].
\end{eqnarray*}\end{proposition}

 Note that each term in the equality
satisfies Poincar\'{e} duality i.e. $f(t)=t^{6g-6}f(t^{-1})$.

\end{document}